\input amstex
\documentstyle{amsppt}
\nologo
\UseAMSsymbols

\magnification=\magstep1

\define\qp{\Bbb Q_p}
\define\1{\sigma_1}
\define\2{\sigma_2}
\define\zp{\Bbb Z_p}

\topmatter

\title Ramification in $C$-extensions of local fields of
characteristic 0\endtitle \rightheadtext{Ramification in $C$-extensions}
\author{Odile Sauzet}\endauthor \address{
University Of Cambridge, Department of Pure Mathematics and Mathematical
Statistics,
Centre for Mathematical Sciences,
Wilberforce Road,
Cambridge CB3 0WB,
UK}
\endaddress 
\thanks Supported by a Marie Curie Fellowship from the European
Union.\endthanks \email{O.Sauzet\@dpmms.cam.ac.uk}\endemail
\abstract We construct explicitly APF extensions of finite extensions of $\qp$
for which the Galois group is not a $p$-adic Lie group and which do not have
any open subgroup with $\zp$-quotient.
 \endabstract 
\endtopmatter

\bigskip Let
$K$ be a finite extension of $\qp$ and $L/K$ be a Galois totally ramified
pro-$p$-extension. If the Galois group of $L/K$ is a $p$-adic Lie group it is
known from Sen ([Sen]) that the sequence of upper ramification breaks of
$L/K$ is unbounded and that the higher ramification groups are open in
$Gal(L/K)$, i.e. in the terminology of [FW], the extension is arithmetically
profinite (APF). 

The more general question concerning the existence of APF extensions of
finite extensions of $\qp$ with a Galois group which is a free pro-$p$-group
remains open (the answer is positive in characteristic $p$, see [Fe]) and we
are interested in obtaining methods to construct  APF extensions with a
non $p$-adic Lie  Galois group.
\medskip
Coates and Greenberg ([CG]) raised the question  whether the notion of
unbounded sequence of upper ramification breaks (extensions deeply ramified)
depends on the existence of an open subgroup of the Galois group having a
quotient isomorphic to $\zp$. Fesenko answered this question negatively in
[Fe] by showing the existence of a finitely generated pro-$p$-group $T$ such
that no open subgroup contains a $\zp$-quotient. He showed that the group $T$
can be realized a the Galois group of a finite extension of $\qp$. 
\bigskip
In this paper we look at  pro-$p$-groups such  that  the lower $p$-series 
define a filtration of open subgroups. We show that this type of group can be
realized as the Galois group of a deeply ramified extension (and arithmetically
profinite). Not only the reason for being deeply ramified is not due to the
presence of a $\zp$-extension but we show that this notion is independent of
the nature of the Galois group, i.e. there are totally ramified extensions with
such a Galois group which are not deeply ramified.

\medskip

We give a method to construct explicitly APF extensions with a
Galois group such as we described above. The method consists in taking
a tower of extension $K_i/K_{i-1}$ such that $K_{i+2}/K_{i-1}$ is a Galois
extension of order $p^3$ with Galois group having two generators $a_i$ and
$a_{i+1}$ and the Galois group of  $K_{i+2}/K_{i+1}$ is generated by the
commutator $[a_{i+1},a_i]$. We show how we can control asymptotically the
ramification in $M=\cup_iK_i/K$. We call the normal closure of such an
extension a $C$-extension (the "$C$" standing for commutators)
  \medskip
After recalling some background on ramification in Section 0, 
 we give results on the structure of the Galois group of the normal
closure of the extension $L/K$ which are relevant to the identification of the
sequence of ramification breaks. Namely, we prove that the Galois
group does not contain any proper closed  normal subgroup. Therefore a
$C$-extension is deeply ramified if and only if it is APF. And we show that
the group does not contain any open subgroup with a $\zp$-quotient.

In Section 2 we construct the tower described above and we show how we can
have it to be APF. In Section 3  we show that the unboundedness of
the ramification breaks does not depend on the structure of the Galois group
because there are $C$-extensions which are not deeply ramified. The last
section shows how one can use $C$-extensions to construct more APF extensions.
Working on problems relative to arithmetically profinite extensions  was
suggested to me by Ivan Fesenko. I am grateful for his support and numerous
comments.

\bigskip

\head 0. Ramification and arithmetically profinite extensions. \endhead
\bigskip\bigskip
We recall some basic notion of ramification theory,  for finite extensions see
[Ser], and for infinite extensions see [FV]. For a general introduction to APF
extensions see [DuFe].
\bigskip
 Let $L/K$ be a totally ramified finite Galois $p$-extension of local
fields  with  finite residue field and with Galois group $G$. Let 
$$i_G(\sigma)=v_L(\frac{\sigma\pi_L}{\pi_L}-1);$$
where $\pi_L$ is a uniformiser of $L$ and $\sigma$ an element in $G$.
\smallskip
The ramification filtration in the lower numbering of $G$ is defined by (see
[Ser1, Chap. VI])  $$G_i=\{\sigma\in G\ ; \ i_G(\sigma)\geq i+1\}.$$
This filtration does not behave well with quotients. 
To define ramification for infinite 
extensions, there is  another 
filtration defined via the Herbrand function $h_{L/K}$ whose explicit
formulation in the case of  totally ramified 
$p$-extensions is as follow :
$$h_{L/K}(x)=\left\{\matrix x\ \ \ \ \text{ if }\ x\leq i\\
(p-1)i+px \ \ \text{ otherwise}\endmatrix\right. \text{ if }[L:K]=p\ ;$$
and for a general finite extension the function is obtained with the equality :
$$h_{L/K}=h_{L/M}\circ h_{M/K}$$
where $M/K$ is a subextension of $L/K$.
\smallskip
Let consider a group $G$ of wild automorphisms of a local field of characteristic
$p$ with finite residue field, we can define the ramification filtration
in a similar way, we write $h_G$ for the Herbrand function.
\bigskip
Then the ramification filtration in the upper numbering is given by
$$G^u=G_{h_{L/K}(u)}.$$
Let $H$ be a normal subgroup of $G$, then $$(G/H)^u=G^u.H/H.$$
Now we can define the upper ramification filtration for infinite p-extensions :

Let $L/K$ be an infinite Galois extension such that 
$\displaystyle{L=\cup_i K_i}$, then $$\displaystyle{G=Gal(L/K)=\varprojlim Gal(K_i/K)}$$ and 
$$G^u=\varprojlim Gal(K_i/K)^u.$$
\bigskip
\proclaim{Definition 0.1} The totally ramified extension $L/K$ is arithmetically profinite (APF) 
if for all $x\geq 0$, the normal subgroup $G^x$ is non trivial and open in $G$
\endproclaim
\bigskip
If the extension $L/K$ is APF, then we can again define the Herbrand function :
$$h_{L/K}(x)=\lim_{i\rightarrow\infty} h_{K_i/K}(x),$$
and the ramification filtration of $G$ in the lower numbering is well defined.
\bigskip
Fontaine and Wintenberger have defined the field of norms $N(L/K)$ for an
APF extension $L/K$  ([FW]) :
$$N(L/K)^\times=\varprojlim K_i^\times\ \ ; \ \ N(L/K)=N(L/K)^\times\cup\{0\}$$
where $\{K_i\}$ is an ordered family of finite subextensions of $L/K$ such 
that $L=\cup K_i$ and the connecting maps are the norm.
If $L/K$ is APF, then the valuation defined in [W2, 2.1.2.], 
addition defined in [W2, 2.1.3.] and the natural multiplication give to 
$N(L/K)$ a structure of local field of characteristic $p$ with residue field
isomorphic to the the one of $K$ and we have the following result on the 
action of $G$ on $N(L/K)$ and the ramification ([W2, 3.3.4.]) :
\bigskip
\proclaim{Proposition 0.2} $G$ acts faithfully on N(L/K) and can be identified 
to a closed subgroup of $Aut(N(L/K))$. Moreover the ramification filtration of 
$G$ can be read in its image in $Aut(N(L/K))$, i.e. if $N_{L/K}(\sigma)$ is 
the image of an element $\sigma\in G$ in $Aut(N(L/K))$, then 
$$G_i=\{\sigma\in G\ ; \ i_{N(L/K)}(N_{L/K}(\sigma))\geq i\}.$$
\endproclaim
\bigskip
The following notations will be used in  the sequel :

- The ramification break $i$ of a finite $p$-extension or an APF $p$-extension 
with Galois group $G$, or a 
group of automorphism $G$ of a local field of characteristic $p$ is defined by
$$i=\sup\{x>0, G=G_x\}.$$
\smallskip
- The upper ramification break of $G$ is $u(G)=h_{L/K}^{-1}(i(G))$.

\bigskip

\head 1. Definitions and basic properties \endhead
\bigskip
\bigskip
A pro-$p$-extension $L/K$ can be defined as a tower of extensions $K_n/K$ such
that $[K_n:K_{n-1}]=p$. If the extension $L/K$ is Galois they 
correspond to a central series of the Galois group. In order to construct
explicitly a pro-$p$-extension with a given Galois group it is convenient to
describe the Galois group in term of generators for the quotients of two
consecutive terms in a central series.
\bigskip
Let recall how a $p$-group can be presented by powers and commutators 
(see [Jo, 13.2]).

Let $G$ be a $p$-group of order $p^n$ and take a central series 
$$G=G_0\subset G_1\supset ...\supset G_n=1$$ such that $|G_i/G_{i-1}|=p$.
Take $a_i$ a non trivial element of $G_{i-1}\setminus G_i$ for all $1\leq i\leq n$. 
There is a power-commutator presentation of $G$ :
$$G=<a_1,...,a_n\ |\ P,C>$$ where the relations are :
$$P\ :\ a_j^p=\prod_{k=j+1}^na_k^{\alpha(j,k)},\ 1\leq j\leq n, \ 1\leq\alpha(j,k) \leq p  ;$$

$$C\ :\ [a_j,a_i]=\prod_{k=j+1}^na_k^{\beta(1,j,k)},\ 1\leq i<j\leq n,\ 1\leq\beta(1,j,k) \leq p  .$$
We call the $a_i$'s the PC-generators (not to confuse with the group 
minimal set of generators).
\bigskip
The notion of PC-presentation can be extended to pro-$p$-groups by taking the
inverse limit of PC-presented $p$-groups
\bigskip
 \proclaim{Definition 1.1} Let $G$ be a 
pro-$p$-group with the PC-presentation $$G=\varprojlim_k <a_1,...,a_k\ |\
P_k,C_k>,$$ The set $\{a_i\}$ of elements of $G$ is called a set of
PC-generators. 
\endproclaim
\bigskip
$C$-extensions are defined by the nature of a PC-presentation of their Galois
group. 
 \bigskip
We need to define the lower central series of a group. Let $\gamma_1(G)=G$
and $\gamma_{i+1}(G)=[\gamma_i(G),G]$. This series
is strictly decreasing, i.e  $\gamma_i(G)/\gamma_{i+1}(G)$ is trivial if and
only if $\gamma_i(G)$ is trivial.
\bigskip
\proclaim{Definition 1.2} The length of an element $a$ in $G$, is the greater
integer $i$ such that $$a\in\gamma_i(G).$$
\endproclaim
\bigskip
A filtration of a pro-$p$-group by open subgroups can be defined as follow
(for a finitely generated pro-$p$-group) :
$$P_0(G)=G,\ P_1(G)=G^p[G,G],\ P_{i+1}(G)=P_i^p(G)[P_i(G),G].$$
It is called the lower $p$-series.

 \bigskip
 \proclaim{Definition 1.3} Let $G$ be a  pro-$p$-group with a
set of PC-generators $$J=\{a_1,a_2...,a_k,...\}.$$ 
\smallskip

A $C$-tower in $G$ is a subset noted $(a_{n_k})$ of $J$, satisfying
$$a_{n_3}=[a_{n_1},a_{n_2}], \text{ for certain } n_1\text{ and } n_2\ ;$$$$
a_{n_4}=[a_{n_3},a_{n_2}],a_{n_{k+1}}=[a_{n_k},a_{n_{k-1}}].$$ \endproclaim
\bigskip

\proclaim{Definition 1.4} a. A $C$-group is an infinite pro-$p$-group
$G$ which
has a set of  PC-generators $I$ containing a $C$-tower $C$, and if $N$ is the
subgroup of $G$ PC-generated by $I\setminus C$, then $$\bigcap_{g\in
G}g^{-1}Ng=1$$
\smallskip
b. A $C$-extension is a Galois pro-$p$-extension whose Galois group is a
$C$-group.
 \endproclaim
\bigskip

Let $L/K$ be a $C$-extension and $G$ its Galois group.
We will note $N(E)/K$ the normal closure of an extension $E/K$ and $M=L^N$ will
be called the extension of $K$ defined by the $C$-tower $C$. Also we will note
$K_i$ the subfields of $M/K$ such that $Gal(K_i/K_{i-1})=<a_i>$ (where we
identify $a_i$ with its image in  $Gal(K_i/K_{i-1})$).

We call an extension of $C_n$-type an extension $E/K$ such that there exists
a $C$-extension $L/K$ such that $L^N=\cup K_i$ as above and $E=K_n$. If $E$ is
an extension of $C_n$-type, we will note $N(E)$ the normal closure of the
extension $E/K$

 \bigskip
\noindent \bf Remark. \rm The existence of  $C$-extensions will be proved in
Section 2. \bigskip
The purpose of the next lemma is to describe the elements in the Galois group
of a $C$-extension which are not in the $C$-tower.
\bigskip

\proclaim{Lemma 1.5} The elements $b_k$ in $Gal(N(K_k)/K_kN(K_{k-1}))$
are such that $$a_k=[b_k,a], \text{ and } b_k\in<[a_{k-1},a_j],2\leq j\leq
k_2>$$ with $a$ an element in $Gal(N(K_k)/K)$ which does not belong to
$$Gal(N(K_k)/K_kN(K_{k-1})).$$
\endproclaim
\bigskip
\demo{Proof} The assertion $a_k=[b_k,a]$ comes from the fact that
$K_kN(K_{k-1})/K$ is not a normal extension, therefore
$$Gal(N(K_k)/K_kN(K_{k-1}))$$ is a subgroup of $Gal(N(K_k)/K)$ which is not
normal. For all $b_k$ in $$Gal(N(K_k)/K_kN(K_{k-1})),$$ there is an element $a$
in $Gal(N(K_k)/K))$ with a non trivial image in $$Gal(N(K_{k-1})K_k/K)$$ such
that $$[b_k,a]\not\in Gal(N(K_k)/K_kN(K_{k-1}))$$
and as $N(K_k)/N(K_{k-1})$ is a normal extension of $K$, we have that 
$$[b_k,a]=a_k.$$
\medskip

We recall how a normal subgroup of a $p$-group is generated (see [Jo]).
Let $G$ be a $p$-group and $H$ be a normal subgroup. Let $I$ be a set of
representative in $G$ of elements in $G/H$ (Schreier transversal). Let $J$ be
a set of generators of $G$. A generating set of $H$ (not necessarily minimal)
is the following $$\{ tx(\bar {tx})^{-1},\ x\in J,\ t\in I\}$$
where $\bar{tx}$ is the representative in $I$ of the image of $tx$ in $G/H$.
\medskip
The element $a_k$ is a generator of the normal subgroup $H_k=Gal(N_k/N_{k-1})$,
but the normal subgroup of $H_k$ not containing $a_k$ define a subextension
of $N_k/K$ which is not normal. So there is an element $b_{k,1}$ such that 
$$a_k=[b_{k,1},a]\in Gal(N_k/K)$$ with $a$ a lift of an element in  
$Gal(N_{k-1}/K)$. 
If the normal subgroup of $H_k$ not containing $b_{k,1}$ define a subextension
of $N_k/K$ which is  normal then $\{a_k,b_{k,1}\}$ define a system of
generators of the normal subgroup $H_k$ of $Gal(N_k/K)$. If not, there exists
an element $b_{k,2}$ such that 
$$b_{k,1}=[b_{k,2},a]\in Gal(N(K_k)/K).$$
This process is finite because $H_k$ is a finite group.
Now we have to show that the set of generators can be taken as in the lemma :
\medskip
Writing explicitly the commutators we obtain
the following  equality, by assuming that whenever a commutator was of length
greater that the length of $a_k$ it was trivial in $Gal(N_{k-1}K_k/K)$ by the
above :

$$a_k=[a_{k-1},a_{k-2}]=[a_{k-1},[a_{k-3},a_{k-4}]]$$$$
=[a_{k-1}^{-1},a_{k-4}^{-1}][a_{k-4}^{-1},[a_{k-1}^{-1},a_{k-4}^{-1}]][a_{k-4},[a_{k-3},a_{k-1}]]
[a_{k-1}^{-1},a_{k-4}][[a_{k-1}^{-1},a_{k-4}],a_{k-3}].$$
The elements  $$[a_{k-3},a_{k-1}]\text{ and }[a_{k-1},a_{k-4}]$$ are non
trivial if and only if $a_k$ is non trivial.
Similarly we obtain that $[a_{k-3},a_{k-1}]$ is non trivial if
$$[a_{k-5},a_{k-1}]\text{ and }[a_{k-1},a_{k-4}]$$ are non trivial.
And therefore $a_k$ is non trivial if for all $2\leq i\leq k-2$,
$[a_{k-1},a_i]$ is non trivial.
\qed
\enddemo
\bigskip

\proclaim{Proposition 1.6} 
Let $L/K$ be a $C$-extension with Galois group $G$. The pro-$p$-group $G$ is
just infinite, i.e. all proper closed normal subgroups of $G$ are
open. \endproclaim
\bigskip
\demo{Proof} Let $H$ be a closed normal subgroup of $G$. If $H$ contains $a_i$
an element in the $C$-tower, then all $a_j$, $j>i$ are in $H$. Indeed, by
normality of the subgroup $[a_i,a_{i-1}]=a_{i+1}$ is in $H$ and the
result is obtained by induction. In Lemma 1.5, we saw that all the
PC-generators of $G$ which are not in the $C$-tower are elements in the group
spanned by $\{[a_{k-1},a_j],2\leq j\leq k_2\}$. If almost all $a_i$'s are is in
$H$, then almost all $[a_{k-1},a_j]$ is also in $H$ by normality of $H$ and
then almost all PC-generator of $G$ is in $H$. Therefore $H$ is an open
subgroup of $G$.

 \medskip

Now assume that no element of the $C$-tower is in $H$. $H$ contains a
PC-generator $b$ of $G$ which is not in the $C$-tower, such that
$a_k=[b,a]$ for a certain integer $k$ and $a$ is an element of $G$. By
normality of $H$ $a_k$ is in $H$ and as above $H$ is open.

  \qed \enddemo
\bigskip
The direct consequence of this result is that it will be rather simple to
check if a $C$-extension is APF, as it is enough to check the unboundedness
of the sequence of upper ramification breaks.

 \bigskip
\proclaim{Corollary 1.7} A $C$-extension
$L/K$ is APF if and only if it is deeply ramified, i.e. if the sequence of
upper ramification breaks of $L/K$ is unbounded.
\endproclaim
\bigskip
\demo{Proof} Let $G=Gal(L/K)$. If $L/K$ is APF then it is deeply ramified.
Assume that sequence of upper ramification breaks of $L/K$ is unbounded
then every
upper ramification group $G^u$ for $u>0$ is non trivial. By Prop. 1.6 all $G^u$
for $u>0$ is also open.
\qed
\enddemo
\bigskip
We show that these extensions will give  examples of APF extensions where
the unboundness of the sequence of upper ramification breaks will not be the
consequence of a $\zp$-quotient in an open subgroup.

\bigskip
 \proclaim{Proposition 1.8} All open subgroup $H$ of $G$ contains no
normal closed  subgroup $N$ such that $$H/N\simeq\zp.$$
\endproclaim
\bigskip
\demo{Proof}Let $H$ be an open subgroup of $G$ and  there exists an element $a$
which is a generator of $H$ and an element of infinite order ([Ze]).  But the
principle of the proof is to show that for all open subgroup $H$ the derived
group $[H,H]$ is of finite index. Therefore there are no quotient of $H$
isomorphic to $\zp$ which is an infinite abelian group.
\medskip
We show that the lower central series are equal to the lower $p$-series, i.e
$$\forall i,\ \gamma_i(G)=P_i(G)=P_{i-1}(G)^p[P_{i-1}(G),G]$$
We have to show that
 $$\gamma_i^p(G)\subset [\gamma_i(G),G].$$
Let $a$ be an element in $G$ of length $\ell$. There is a commutator $[b,c]$ of
length $\ell$ such that $a=[b,c]a'$.
$$a^p=[b,c]^p{a'}^p...$$
the $...$ stand for commutators of higher length.
$$[b,c]^p=[b^p,c][b,[b,c]]^{p-1}...$$
as we have proved in Lemma 1.5 that there are no subgroup of any finite
quotient of $G$ generated by an element $a_i^p$ such that the length of
$a_i^p$ is the length of $a_i$. So either $b^p=1$ or we can decompose again in
the same way, and after a finite number of time we obtain that the length of
$a^p$ is strictly greater than $\ell$.

 \medskip 
$G$ being a finitely generated pro-$p$-group such that $G^p$ is included
in the derived group $[G,G]$, this last subgroup is of finite index in $G$. Let
$H$ be an open subgroup of $G$, there exists an integer $i$ such that 
$$\gamma_{i}(G)\subset H.$$
Therefore 
$$[\gamma_{i}(G),\gamma_{i}(G)]\subset [H,H].$$
To show that the derived group $[H,H]$ is of finite index in $H$ it is enough
to show that for all $i$,  $[\gamma_{i}(G),\gamma_{i}(G)]$ is of finite index
in $G$. 

There is an inclusion  $$[\gamma_{i}(G),\gamma_{i}(G)]\subset\gamma_{2i}(G).$$
The quotient $$\gamma_{2i}(G)/[\gamma_{i}(G),\gamma_{i}(G)]$$ contains only
elements in 
$$[G/\gamma_2(G),\gamma_{2i-1}(G)/\gamma_{2i}(G)],[\gamma_2(G)/\gamma_3(G),\gamma_{2i-2}(G)/\gamma_{2i-1}(G)],...,$$$$
[\gamma_{i-1}(G)/\gamma_i(G),\gamma_{i+1}(G)/\gamma_{i+2}(G)].$$ 
All these groups being finite, the quotient 
$$\gamma_{2i}(G)/[\gamma_{i}(G),\gamma_{i}(G)]$$ is finite. The subgroup
$\gamma_{2i}(G)$ is of finite index in $G$, and so is
$[\gamma_{i}(G),\gamma_{i}(G)]$. For every open subgroup $H$ of $G$, the
derived group $[H,H]$  contains a subgroup which is open in $G$ and we have
proved the result.  \qed \enddemo

\bigskip
The last result of this section shows that if a $C$-tower is APF, this is not
due to the result of Sen on $p$-adic Lie groups.

\bigskip
\proclaim{Proposition 1.9} The Galois group of a $C$-extension 
is not a $p$-adic analytic group.
\endproclaim
\bigskip
\demo{Proof} It is enough to show that there are open subgroups with an 
arbitrary large number of generators, i.e we show that the pro-$p$-group $G$ is
not of finite rank (see [DDMS] Prop. 3.1 and Cor. 8.34). Consider the subgroup
$H$  of $G$ containing 
$$a_2,a_3,...,a_{2k-1},a_{2k+1},a_{2k+2},...$$ ie. for $i<2k+1$ even, $a_i$ is
not in the subgroup and for any $i\geq 2k+1$ $a_i$ is in the subgroup.
Then the elements $$a_2,a_3,...,a_{2k-1},a_{2k+1}$$ are topological
generators of the open subgroup $H$. Indeed, none of this elements are in the
Frattini subgroup $\phi (H)$ or are equal to another one modulo an element in
$\phi (H)$ : $$a_{2i+1}=[a_{2i},a_{2i-1}], a_{2i}\not\in H, a_{2i-1}\in H.$$

 \qed
\enddemo

\bigskip

\bigskip
\head 2. There exist APF $C$-extensions \endhead
\bigskip
\bigskip
In this section we assume that $K$ is a finite extension of $\qp$ of degree
at least 3. We construct a totally ramified $C$-extension and show that there
exist $C$-extensions of $K$ which are APF. \bigskip

 \proclaim{Proposition 2.1} Let $K$ be a finite  extension 
of $\qp$ of index of ramification $e(K)$. For all  $n>0$ there exists an
extension of $C_n$-type with maximal upper ramification break satisfying 
$$u_n>ne(K)-v_n\text{ where } \lim_{n\rightarrow\infty}v_n<\infty.$$
\endproclaim
\bigskip
First we prove a lemma.
\bigskip
\proclaim{Lemma 2.2}
i. Let $K_1/K$ and $E/K_1$ be cyclic extensions of order $p$ such that  $E$ is
not normal over $K$. Then there exists  a cyclic  extension $K_2/K$ of order
$p$  such that $K_2E/K$ is a normal extension of order $p^3$ with Galois group
having the presentation : $$Gal(K_2E/K)=<a_1,a_2\ |\
a_1^p=a_2^p=[a_1,a_2]^p=1, [a_1,[a_1,a_2]]= [a_2,[a_1,a_2]]=1>$$

ii.  Let $K_1/K$ be a cyclic extension of order $p$ and $E/K_1$ an extension 
of $C_n$-type. We have one of the following cases :

a.  there exist
 a cyclic 
extension $K_2/K$ of order $p$ such that $K_2E/K$ is a normal extension
of $C_{n+2}$-type ;

b.  there exists
 two cyclic 
extensions $K_2/K$ and $K_3/K$ of order $p$ such that
$K_2K_3N(E)/K$ is a normal  extension and if we note $b$ and $c$ a lift of the
generators $Gal(K_2/K)$ and $Gal(K_3/K)$ in   $$Gal(K_2K_3N(E)/K)$$ then
we have the relations in $Gal(K_2K_3N(E)/K)$

$$[a,b]=a_1\alpha_1\ \ ;\ \ [a,c]=a_2\alpha_2$$
where $(a_1\alpha_1,a_2\alpha_2)$ define a minimal system of generators of 
$$Gal(K_2K_3N(E)/K_1K_2K_3)$$ and $a$ is a lift of a generator of
$Gal(K_1/K)$.
 \endproclaim
\bigskip
\demo{Proof} i. Let take $E(K)$ the normal closure of $E/K$ and consider a 
generator $b$ of $Gal(E(K)/K)$ which has a trivial image in $Gal(K_1/K)$. We 
note $a$ (resp. $c$) a lift of a generator of $Gal(K_1/K)$ (resp. $Gal(E/K_1)$)
in $Gal(E(K)/K)$.
Then the subgroup generated by $b$ is not a normal subgroup. It means that one
can choose the generators such that  there is one of the relations which is 
$$[a,b]=c\ \text{ or }\ [b,c]=a.$$
We note that the  subgroup generated by $b$ and $c$ is normal so the 
second relation is not possible and we get the result.
\bigskip
ii. 
We obtain two extensions $K_2/K$ and $K_3/K$ by considering
two cyclic independent subextensions of $E/K_1$ as in case i.
Let call $b$ and $c$ the two generators obtained. There are two 
possibilities : either $K_2$ and $K_3$ are linearly independent in which case
we are in case b. or $c=b^i$ $i\not=1$ and we are in case a, indeed we have
that
$a_1=[a,b]$ and $$a_2=[a,b^i]=[a,b]^i[[a,b],b]$$ and therefore
$$[a_2,a_1]=[[[a,b],b],[a,b]].$$

\qed
\enddemo

\demo{Proof of Proposition 2.1} 
Let show first that for all $K$ non trivial extension of $\qp$, there exists an
extension $E_3/K$ of  $C_3$-type satisfying the condition on the ramification.
Let take an extension $K_1/K$ with ramification break $i_1(3)$ and an extension
$E/K_1$ with  ramification break $i$ all of them cyclic of order $p$ such
that $i>i_1(3)$ and  we put 
$$u_3=i_1(3)\frac{p-1}{p}+\frac{i}{p}.$$
Moreover we want $E$ not to be normal over $K$. For this it enough to take
$i$ such that $p$ does not divide $i-i_1$.
By Lemma 2.2, there exists a normal extension $E_3/K$ of $C_3$ type for
which the bigger upper ramification break is  $u_3$.
\bigskip
Let $K_1/K$ be a cyclic extension of order $p$ with ramification breaks
$i_1(5)$. We showed that there exists an extension $E_3/K_1$ with bigger upper
ramification break $u_3$. We put
$$u_5=i_1(5)\frac{p-1}{p}+\frac{u_3}{p}.$$
By Lemma 2.2
we have two possibilities either $c=b^i$ and we have the result or
$Gal(E(K)/K)$ is a $p$-group with exactly three generators.

\medskip 
We show for all $n$ that $Gal(E(K)/K)$ contains a  subgroup $N$
such that the subextension $E(K)^N/K)$ is the subextension of an
extension of $C_n$-type.

Take the normal subgroup $N$ of $Gal(E(K)/K)$ generated by $c[a,b]$. The
group $N$ also contains $$a^{-1}c[a,b]a=c[c,a][a,b][[a,b],a]\text{ and}$$
$$b^{-1}c[a,b]b=c[c,b][a,b][[a,b],b]$$
which will define relations in $Gal(E(K)^N/K)$.

The map $$\matrix Gal(E(K)/K)&\rightarrow &Gal(E(K)/K)/N\\
                    a&\rightarrow&\bar a\\
                     b&\rightarrow &\bar b\\
                    c&\rightarrow&[\bar b,\bar a]\\
                    a_1&\rightarrow&[\bar a,\bar b]\\
                     a_2 &\rightarrow&[\bar a,[\bar a,\bar b]]\\
             b^{-1}c[a,b]b&\rightarrow&1\\
a^{-1}c[a,b]a&\rightarrow&1=[\bar b,\bar a][[\bar b,\bar a],\bar a][\bar
a,\bar b][[\bar a,\bar b],\bar a] \\
...&&\endmatrix$$ 
shows that $E(K)^N/K$ is the subextension of an extension of $C_n$-type
because  the relation

$$[\bar b,\bar a][[\bar b,\bar a],\bar a][\bar
a,\bar b][[\bar a,\bar b],\bar a]=1$$
is trivial  :
$$a_3^{-1}[a_3^{-1},a_2]a_3a_4=a_2^{-1}a_3^{-1}a_2a_3a_4=a_2^{-1}a_3^{-1}a_3a_2a_4^{-1}a_4=1$$
Let consider the extension of $C_n$-type which has $E(K)^N$ as a subextension.
We proceed with the construction. The bigger ramification break in that
extension is $u_5$ defined above.
\medskip
Now we proved the result for $n=5$ and we assume we can show it for any odd
integer $n$ and show that then for all finite extension $K$ of $\qp$
there exists an extension of $C_{n+2}$-type with bigger upper ramification
break $u_{n+2}$.

Let $K_1/K$ be a cyclic extension of order $p$ with ramification break
$i_1({n+2})$. By assumption there exists an extension of $K_1$ of $C_n$-type
with biggest upper ramification break $u_n$
then put
$$u_{n+2}=i_1({n+2})\frac{p-1}{p}+\frac{u_n}{p}.$$

By Lemma 2.2, the situation is similar to the case $n=5$. Either there
exists an extension $K_2/K$ such that $E_nK_2/K$ is of $C_{n+2}$ type or we
have to take a quotient of $Gal(E_n(K)/K)$ by the same normal subgroup as in
the previous case. 
\bigskip
It remains to show that  the ramification breaks
$i_1(3),i_1(5),...,i_1(2n+1),...$ can be taken such that upper ramification
breaks $u_3,..,u_n,...$ satisfy the condition in the proposition, i.e there
exists a convergent sequence $v_n$ such that $$u_n>ne(K)-v_n.$$
For $n=3$ take any $i_1(3)$ and $i$.

Then for $n=5$, we take $$i_1(5)=\frac{pe(K)}{p-1}-\varepsilon_5$$ where
$\varepsilon_5$ is an integer, such that 
$$u_5=e-\frac{p-1}{p}\varepsilon_5+\frac{u_3}{p}.$$

For $n=7$, we take $$i_1(7)=\frac{pe(K)}{p-1}-\varepsilon_7$$
such that 
$$u_7=e-\frac{p-1}{p}\varepsilon_7+\frac{u_5}{p}=
e-\frac{p-1}{p}\varepsilon_7+\frac{pe(K)}{p}-\frac{p-1}{p^2}\varepsilon_5+\frac{u_3}{p^2}$$$$=
2e-[\frac{p-1}{p}(\varepsilon_7+\frac{\varepsilon_5}{p})-\frac{u_3}{p^2}].$$

Note that for $n+1$ we have that the value $e(K_1)=pe(K)$
\medskip
By taking for all odd integer $k$, $$i_k=\frac{pe(K)}{p-1}-\varepsilon_k$$ we
obtain
$$u_{2n+1}=(n-1)e-[\frac{p-1}{p}(\varepsilon_{2n+1}+\frac{\varepsilon_{2n-1}}{p}+...
+\frac{\varepsilon_5}{p^{n-1}})-\frac{u_3}{p^{n-1}}].$$
Take all the $\varepsilon_{k}$ in a bounded interval, then 
$$v_n=[\frac{p-1}{p}(\varepsilon_{2n+1}+\frac{\varepsilon_{2n-1}}{p}+...
+\frac{\varepsilon_5}{p^{n-1}})-\frac{u_3}{p^{n-1}}].$$
is a convergent sequence.

\qed
 \enddemo
\bigskip
\proclaim{Corollary 2.3} There exist APF $C$-extensions.
\endproclaim
\bigskip
\demo{Proof}
Let $L/K$ be a $C$-extension and note $\{a_i\}$ a $C$-tower defining  its
Galois group. We assume that this extension is not APF. Take a family $E_n$ of
extensions of $K$ linearly independent of $L/K$ such that for all $n$ $E_n/K$
is an extension of $C_n$-type and its greater upper ramification break $\tilde
u_n$ satisfies $$\tilde u_n\geq ke\text{ where } n=2k \text{ or } n=2k-1.$$
\medskip
This family exists by Proposition 2.1 and one can check that the extensions
taken are linearly independent of $L/K$ by checking that the compositum of the
$L^{ab}/K$ and the subextensions $E_{n,1}/K$ of $E_n$ of order $p^2$ which is 
not normal remains not normal. 
 \medskip
We will note $\{\tilde a_{i,n}\}$ the $C$-tower
defining the extension $E_n/K$. 
 \medskip Let $u_n$ be the sequence
of upper ramification breaks in $L/K$. This sequence is convergent as we
assumed that  the extension is not APF. We  modify the
extension $L/K$ in order to obtain a $C$-extension $\tilde L/K$ which is APF.
Consider the compositum $M$ of $L$ with  all the finite extensions $
N(E_n)/K$  and take $\tilde L$ the  fixed field of $M$ by the normal subgroup
generated by $a_1\tilde a_{1,n}^{-1}$ and $a_2\tilde a_{2,n}^{-1}$ for all $n$.
We show that this is a well defined $C$-extension and that it is APF, i.e. that
the sequence of upper ramification breaks in $\tilde L/K$ is unbounded (see
Cor. 1.7). 
\bigskip
The Galois group $Gal(\tilde L/K)$ is defined by the surjective map
$$\matrix Gal(M/K)&\rightarrow& Gal(\tilde L/K)\\
             a_1&\rightarrow& \bar a_1\\
             a_2 &\rightarrow& \bar a_2\\
             a_{1,n}&\rightarrow& \bar a_1\ &\forall n>2\\
               a_{2,n}&\rightarrow& \bar a_2&\forall n>2\\
                 a_k   &\rightarrow& \bar a_k( \bar a_1,  \bar a_2)&\forall n>2\\
               a_{k,n}&\rightarrow&   \bar a_k( \bar a_1,  \bar a_2)&\forall n>2
     \endmatrix$$
where $\bar a_k( \bar a_1,  \bar a_2)$ is the $k^{th}$ element is the
$C$-tower generated by $\bar a_1$ end  $\bar a_2$. This a group homomorphism.
So $\tilde L/K$ is a $C$-extension.
\medskip
The upper ramification breaks $\bar u_k$ in $\tilde L/K$, for $k$ large
enough satisfy:
 $$\bar u_k\geq Sup\{u_{L/K}(a_k),u_{E_n/K}(a_{k,n}),\forall n\geq
k\}.$$ 
Indeed, this is not the case if for all $n\geq k$
$$u_{L/K}(a_k)=u_{E_n/K}(a_{k,n})\text { and
}u_{M/K}(a_ka_{k,n})>u_{L/K}(a_k).$$
If $L/K$ is not APF, for $k$ big enough
$$u_{L/K}(a_k)<u_{E_k/K}(a_{k,k})$$ as
the first sequence is convergent and
$$u_{E_k/K}(a_{k,k})>([\frac{k-1}{2}]-1)e(K).$$
With the hypothesis made on the greater upper ramification break in the
extension $E_n/K$ we obtain that for $n$ big enough
$$\bar u_n\geq ([\frac{n-1}{2}]-1)e(K).$$
and therefore the extension $\tilde L/K$ is APF.
\qed
\enddemo
\bigskip
\bigskip
\head 3. There exist non arithmetically profinite $C$-extensions  \endhead
\bigskip
\bigskip
We want to show that unlike  extensions whose Galois group is a $p$-adic Lie
group, the structure of the Galois group of $C$-extensions does not implies
that the extensions is APF. Indeed, we show how to construct a $C$-extension
which is not APF.

 \bigskip
\proclaim{Proposition 3.1} There exist $C$-extensions which are not deeply
ramified
. \endproclaim
\bigskip
We construct a $C$-extension $L/K$ such that  $L^N/K$ has bounded ramification.
\medskip
 Let $L/K$ be a $C$-extension such that
$L^N=\cup K_n$ as defined in Section 1.
Consider the (not normal) closed subgroup $H$ of $G=Gal(L/K)$ such that only
the elements of the $C$-tower $a_{2i+1},i\geq 0$ are not in $H$ 
and take $M'=L^H$  the fixed field of $L$ by $H$. Then
$M'/K$ is an infinite extension with $M'=\cup M_n$ such that for all $n$
$M_n/M_{n-1}$ is a cyclic extension of order $p$ and $M_n$ is not a normal
extension of $M_{n-2}$. Moreover we assume that for all $n$ we have that the
ramification breaks satisfy : $$p\nmid (i(M_n/M_{n-1})-i(M_{n-1}/M_{n-2})).$$
This is possible to obtain by the methods used in Section 2.

 \bigskip
We are going to construct a family of extensions $E_n/K$ of $C_n$-type for all
$n>1$ such that there exists a subextension of $L\cup(\cup N(E_n))$  which is a
non APF $C$-extension.
\bigskip
First we show that we can construct a tower of extensions $E'_i/K$ such that
for all $i$  the extension $$M_iE'_{i-1}E'_iM_{i-1}/M_{i-1}E'_{i-1}$$
contains a cyclic subextension of order $p$ with a given ramification break
$$i_n<i(M_iE'_{i-1}/M_{i-1}E'_{i-1}).$$

\bigskip
We show a general lemma on the existence of Galois extension of order $p^2$
with given ramification breaks : 
\bigskip
\proclaim{Lemma 3.2}
 Let $L/K$ be a cyclic extension of order $p$ with 
ramification break $i$. Let $j$ and $s$ be two positive integers not 
divisible by $p$ such that if $s\leq i$ then
$s\leq j$ and $i+p^{-1}(s-i)\leq j$ otherwise ; and $\displaystyle{j\leq\frac{pe}{p-1}}$ 
where $e$ is the absolute index of ramification of $K$. If $i\not=j$, take 
$s=h_{L/K}(j)$. There exists a cyclic 
extension $L'/K$ of 
order $p$ with 
ramification break $j$ such that the ramification break of the cyclic 
extension $LL'/L$ is $s$.
\endproclaim
\bigskip
\demo{Proof}Suppose that $i\not=j$.

There exists a cyclic extension $K_1/K$ of order $p$ totally ramified with ramification 
break $j$ if and only if $\displaystyle{j\leq\frac{ep}{p-1}}$ ([FV, Prop. III2.5]). 
As $i\not=j$, the ramification in $LL'/L$ is $s=h_{L/K}(j)$.

\bigskip

Suppose that $i=j$. Let $\lambda$ be a root of the polynomial
 $X^p-X-\alpha$ with $v_K(\alpha)=-i$ and $L=K(\lambda)$ and $<\tau>=Gal(L/K)$.
We show that there exists an element $\beta$ in $K$, $v_K(\beta)=-i$  and
$\mu$  the root of the polynomial
$X^p-X-\beta$ such that
there is an element $A$ in $K(\lambda,\mu)$ stable by $\tau\sigma$ where 
$\sigma$ generates $Gal(K(\mu)/K)$ satisfying the conditions :

. $v_{K(\lambda,\mu)^{\tau\sigma}}(A)=-s$

. $v_{K(\lambda,\mu)^{\tau\sigma}}(\sigma(A)-A-1)>0$.

It is enough to show that $i(Gal(K(\lambda,\mu)^{\tau\sigma}/K))=-s$ and 
therefore
$$i(Gal(K(\lambda,\mu)/L))=-s.$$
\medskip
If $p=2$ take $$A=\mu\tau(\lambda)+\lambda\sigma(\mu),$$ if $p>2$ take
$$A=\sum_{0\leq i<j\leq p-2}\tau^i(\lambda)\sigma^j(\mu)-
\tau^j(\lambda)\sigma^i(\mu).$$
\medskip
First we show that $\tau\sigma(A)=A$.
$$\tau\sigma(A)=\sum_{0\leq i<j\leq p-2}\tau^{i+1}(\lambda)\sigma^{j+1}(\mu)-
\tau^{j+1}(\lambda)\sigma^{i+1}(\mu).$$
We conclude by the fact that $Tr_{K(\mu)/K}(\mu)=Tr_{K(\lambda)/K}(\lambda)=0$ 
which allows to replace $\tau^{p-2+1}(\lambda)$ by 
$\displaystyle{-(\sum_{i=0}^{p-2}\tau^i(\lambda))}$ and 
$\sigma^{p-2+1}(\mu)$ by $\displaystyle{-(\sum_{i=0}^{p-2}\sigma^i(\tau))}$.

For $p=2$, the result is clear.
\bigskip

We show that the second condition is always valid.
We have that $$v_{K(\lambda)}(\tau^j(\lambda)-\tau^i(\lambda)-1)>0$$
and 
$$v_{K(\mu)}(\sigma^j(\mu)-\sigma^i(\mu)-1)>0$$
for $j>i$. For $p>2$
$$\sigma(A)-A=\sum_{0\leq i<j\leq p-2}\tau^i(\lambda)(\sigma^{j+1}(\mu)
-\sigma^j(\mu))-\tau^j(\lambda)(\sigma^{i+1}(\mu)-\sigma^i(\mu))$$
$$\sigma(A)-A=\sum_{0\leq i<j\leq p-2}(\tau^i(\lambda)-\tau^j(\lambda)+...)$$
where the dots mean terms of higher valuation.
The above sum is then equal to
$$\sigma(A)-A=1-p+ B$$ where $B$ is of positive valuation, which gives the 
result.
\bigskip
In the case $p=2$ one has
$$\sigma(A)-A=(\sigma(\mu)-\mu)(\tau(\lambda)-\lambda),$$
so we can conclude in a similar way.
\bigskip

Let $M=K(\lambda,\tau)^{\tau\sigma}$. By a similar method as above we obtain : 
$$v_M(A)=v_M(\sum_{0\leq i<j\leq p-2}\tau^i(\lambda)(\sigma^j(\mu)-
(\sigma^i(\mu))-\sigma^i(\mu)(\tau^j(\lambda)-\tau^i(\lambda)))$$
$$=v_M(\sum_{0\leq i<j\leq p-2}(\tau^i(\lambda)-\sigma^i(\mu))).$$
From the equality of the valuations
$$v_M(\tau^i(\lambda)-\sigma^i(\mu))=v_M(\tau^i(\lambda)^p-\sigma^i(\mu)^p-\alpha+\beta)$$
we deduce that $v_M(\tau^i(\lambda)-\sigma^i(\mu))=-s$ if
$$v_K(\alpha-\beta)=-s\tag{*}$$ and if this last equality is valid then
$v_M(A)=v_K((p-2)(\alpha-\beta))=-s$. So take $\beta=\alpha + \pi_K^{-s}$ and
we have proved the existence of the extension $K(\mu)$ with the required
properties. \medskip For $p=2$ we obtain the condition $v_M(A)=-s$ if
$v_K(\alpha +\beta)=-s$. \qed
\enddemo

\bigskip
The lemma gives  the existence of an extension $F/M_{i-1}E'_{i-1}$
satisfying the property, but we have to make sure that this correspond to 
 an extension of $E'_{i-1}$. For this take an
element $\beta$ in $E'_{i-1}$ such that (see (*) in the
proof of the lemma) $$v_{M_{i-1}E'_{i-1}}(\alpha-\beta)=-s.$$
It remains to check that the extension obtained is not normal over $E'_{i-1}$.
This is clear with the assumption on the lower ramification breaks of the
extensions $M_n/M_{n-1}$ made at the beginning of this section.

 \bigskip
We have a tower of extensions $E_i'$ such that $E'_{i+1}/E'_{i-1}$ is
not normal. We construct the extensions $E_n/K$ in a similar way that we
constructed the extensions $E_n/K$ in the previous section.
\medskip
The extension $E_3/K$ is obtained as in Section 3 by taking the normal
closure of $E'_2/K$. Then we obtain $E_5$, we start by taking the normal
closure of $E'_3/E_1'$ to obtain a Galois extension $F_5/E_1'$ which is  not
normal over $K$ and then take the the normal closure of this extension. We get
$E_5/K$ as in Section 2.
\medskip
We can repeat this process for all $n$ and obtain the required family of
extensions. Now we look at the ramification in the subextension $L'$ of 
$L\cup(\cup N(E_n))$ which is the fixed field of $L\cup(\cup
N(E_n))$ by the normal subgroup generated by $a_1\tilde a_{1,n}^{-1}$ and
$a_2\tilde a_{2,n}^{-1}$ for all $n$ keeping the notation of Section 2. We saw
that it was a $C$-extension. We note that for all $i,n$
and $m$, $\tilde a_{2i,n}=\tilde a_{2i,m}$. In the construction of the tower
of $E'_n$, we can choose the ramification breaks $i_n$ to be small enough so
that the sequence of upper ramification breaks in the extension $L'/K$
converges.  \bigskip

\head 4. Construction of some APF extensions whose Galois group is not a
$p$-adic Lie group \endhead
\bigskip
\bigskip
First we discuss deeply ramified extensions. We know that all totally
ramified extensions whose Galois group contains an open subgroup with a
$\zp$-quotient are deeply ramified. In this paper we have shown that :
\bigskip
\proclaim{Proposition 4.1} i. All pro-$p$-groups that contains a $C$-tower can
be realized as the Galois group of a deeply ramified extension. 

ii. The fact that the
Galois group of a totally ramified extension contains a $C$-tower is not
enough to say that the extension is deeply ramified.
\endproclaim
\bigskip

Now we look at arithmetically profinite extensions.
\smallskip

The result obtained in Section 2 can be extended to construct more APF Galois
pro-$p$-extensions \bigskip
\proclaim{Proposition 4.2} Let $E/K$ be a totally ramified Galois
pro-$p$-extension such that 

. $Gal(E/K)$ is a pro-$p$-group over two
generators ;

. $E$ is the compositum of a finite number of extensions $E_i/K$
each of them containing a subextension $L_i/K$ such that $E_i/L_i$ is a
$C$-extension. 
\medskip
Then for all subextension $L/K$ of $E/K$ which is infinite,
$Gal(N(L)/K)$ contains a  subgroup which is generated by a $C$-tower (i.e
there is a finite extension $K'/K$ such that $L/K'$ is a $C$-extension.
\endproclaim
\bigskip
\demo{Proof} Let call $$\{a_{k,i},i>0\}\text{ for }1\leq k\leq n$$ the
$C$-towers generating a lift of $Gal(E_i/L_i)$ in $G=Gal(E/K)$.
\smallskip
Let $H$ be a normal subgroup of $G$ of infinite index. We show that at least
one $C$-tower is not contained in $H$. Suppose that $$a_{k,i}\in H$$ for a
certain pair $(k,i)$ and $$a_{k,i-1}\not\in H.$$
then as $H$ is normal we have
$$a_{k,i-1}^{-1}a_{k,i}a_{k,i-1}=a_{k,i}a_{k,i+1}\in H$$ 
and therefore $$a_{k,i+1}\in H.$$
So as there are only a finite number of $C$-towers and nothing else there is a
$C$-tower which has a non trivial image in $G/H$.
\medskip
We show that there is an open normal subgroup of $G/H$ which is generated by a
$C$-tower. Let $\{\bar a_i\}$ be the image of a $C$-tower in $G/H$. As $G$
is the union of a finite number of $C$-towers, there are a finite number of
$\bar a_i$ such that there exists an element $\bar b$ in $G/H$ which is not in
the subgroup generated by $\{\bar a_i\}$ which does not commute with $\bar
a_i$. So there exists $i_0$ such that the subgroup generated by the $C$-tower
$\{\bar a_i,\ i>i_0\}$ is an open normal subgroup of $G/H$.

\qed
\enddemo
\bigskip
The sequence of upper ramification breaks in a quotient depends on the
relations that define this quotient  :
 \bigskip
\proclaim{Proposition 4.3} Let $E/K$ be as in Prop. 4.2. Let $H$ be a
closed normal subgroup of infinite index. Let $\{\bar a_i\}$ be a $C$-tower in
$G/H$. The upper ramification breaks of the elements in $\{\bar a_i\}$ for
$i$ big enough are given by :
$$u_{G/H}(\bar a_i)=Max\{u_G(a_{k,i}),\varphi(a_{k,i})=\bar a_i\}.$$
\endproclaim
\bigskip

 \bigskip
\demo{Proof} To see this we have to show that if there exists a $C$-tower
$\{a_{k,i}\}$ which is send to  $\{\bar a_i\}$ in $G/H$, then either no other
$C$-tower is mapped to this one or if there is then 
$a_{l,i}...$ is mapped to $\bar a_i$ where the dots stand for elements
$a_{l,j}$, $j>i$.
 \medskip
This comes from the fact that if $$a_{k,i}\prod_{l\in I,j\in J}a_{l,j}\in H$$
then the towers $\{a_{l,i}\}$ have a non trivial image in $G/H$ if an elements
that does not commutes with $a_{l,j}$ is mapped to $\bar a_i$.
\medskip
We can show that this element can only be $a_{l,j-1}...$. And by going down
until there are no more elements in $G/H$ which commute with one of the
elements, we see that we must have $i=j$

\qed
\enddemo
\bigskip
We have the following result on the ramification in the extension
$L/K$ :
  \bigskip
\proclaim{Corollary 4.4} If for all $C$-towers in $G$ we have that
$$u_G(a_{k,i}a_{l,i})=Inf\{u_G(a_{k,i}),u_G(a_{l,i})\}$$ and all the
$C$-extensions $E_i/L_i$ are APF. Then $E/K$ is APF.
\endproclaim

A consequence of this result is that we can construct an APF extension with a
Galois group as in Proposition 4.3
\bigskip
\proclaim{Proposition 4.5} There exists extension $E/K$ as in Proposition 4.2
which are arithmetically profinite.

 \endproclaim

\bigskip
\demo{Proof}
To construct extensions $E/K$ which are APF we use the same method as in
Section 2. By Corollary 4.4, it is enough to take all the $C$-towers
 having distinct upper ramification breaks. 

Take an extension $E/K$ as in Prop. 4.2. We can modify it so that all the
$C$-extensions $E_i/L_i$ are APF as done in Section 2. In order to have the
condition of Cor. 4.4 we have to modify them so that they have all distinct
upper ramification breaks sequences which can be done because there is a finite
number of $C$-extensions to consider.
\qed
\enddemo
\bigskip
A similar result can be obtained for the compositum of infinitely many
$C$-extensions over the same two generators if we restrict to such group
in which every  elements in a $C$-tower commutes with almost all elements
in $G$. But even in that case there are infinite Galois subextensions whose
Galois group does not contain any $C$-tower. Namely one can take the normal
subgroup that contains all the $a_{k,i}$ for $i$ greater that a given $i_0$,
for all $k$. We can show that there no other Galois subextensions than those
and the one which we looked at the beginning of that section but we have to
make sure that such subextensions are also APF.

We discuss briefly how to generalize the construction made in Section 2, in
order to obtain APF extensions of this type.
\bigskip
Let $K_1/K$ be a cyclic extension of order $p$ and let $K_2/K_1,...,K_p/K_1$
be $p-1$ distinct cyclic extensions of $K_1$ with different ramification
breaks such none of them are normal over $K$. Then the result obtained in
Lemma 2.1 extend and therefore there is a quotient of the normal closure of
the extension $K_2K_3...K_{p}/K$ whose Galois group is PC-generated by
$a_1,a_2,a_3,...a_{p+1}$ with $a_3=[a_2,a_1]$ and $a_i=[a_{i-1},a_1]$ for
$i>3$. 
\bigskip
Now we have to check that, in this type of construction, we can take the
ramification to be hight enough so that  the whole extension is APF. By
considering the number of generators in an open subgroup of a free
pro-$p$-group with two generators, the bound on the upper ramification and the
fact that we want the upper ramification breaks to be unbounded, we can't
simplify our problem by taking all ramification distinct because there would
not be enough integer available. But as seen in Section 3, we can control the
extensions so that the ramification breaks would not collapse in certain
quotients.


\bigskip
\bigskip












\Refs\widestnumber\key{DDMS}
\ref \key CG \by J. Coates \& R. Greenberg Kummer theory for abelian varieties
over local fields \jour Invent. Math. \vol 124 \yr 1996 \pages 129--174 \endref
 \ref \key DDMS \by J.D. Dixon, M.P.F. du Sautoy, A. Mann, D.
Segal \paper Analytic  pro-$p$-groups,   \paperinfo Cambridge studies in
advanced mathematics \vol 61 \jour 2nd ed, Cambridge Univ. Press \yr 1999
\endref

\ref \key DuFe \by M.P.F. du Sautoy \& I.B. Fesenko \paper Where the wild 
things are : ramification groups and the Nottingham group  \jour   New
horizons in pro-$p$-groups,du Sautoy and A.
Shalev (eds.) , Progress in Math. \vol 184  \yr 200\pages 287--328 \endref

\ref \key Fe \by I.B. Fesenko \paper On just infinite pro-$p$-groups and 
arithmetically profinite extensions of local field\jour J. reine angew. Math.
\vol 517\yr 1999\pages 61--80 \endref
\ref \key FV \by I. Fesenko \& S. Vostokov \paper Local Fields and their 
Extensions, a Constructive Approach \jour Transl. of Mathematical Monographs 
\vol 121\endref
\ref \key FW \by J.-M. Fontaine \& J.-P. Wintenberger \paper Corps de normes de 
certaines extensions algebriques de corps locaux \jour Comptes rendus \vol
288  s\'erie A \yr 1979 \pages 441--444 \endref 

\ref\key Jo \by D.L. Johnson \paper Presentation of groups \jour LMS Student Texts \vol 15 \yr 1990 \endref

\ref \key Sen\by S. Sen \paper Ramification in $p$-adic Lie extensions \jour
Invent. Math. \vol 90 \yr 1972 \pages 44--50 \endref

\ref \key Ser \by J.-P. Serre \paper Corps Locaux \jour Hermann \yr 1968
\endref


\ref \key Ze \by E.I. Zelmanov \paper On periodic compact groups \jour Israel
J. Math. \vol 77 \yr 1992 \pages 83--95 \endref
\endRefs
\enddocument